\documentclass[a4paper, 12pt]{article}

\usepackage{amssymb}
\usepackage[definethebibliography]{easybib}
\usepackage{amsthm}
\usepackage{amsmath}
\usepackage{mathrsfs}
\usepackage{color}

\begin{document}

\newtheorem{thm}{Theorem}[section]
\newtheorem{lem}[thm]{Lemma}
\newtheorem{rem}[thm]{Remark}
\newtheorem{cor}[thm]{Corollary}
\newtheorem{prop}[thm]{Proposition}


\renewcommand{\theequation}{\arabic{section}.\arabic{equation}}
\def\proof{\noindent{\it Proof.\ }}
\newcommand{\question}[1]{{\par\medskip\hrule\medskip\noindent {\small\sf #1} \medskip\hrule\medskip\par}}

\def\qed{\hfill $\square$}

\newcommand{\R}{\mathbb R}
\newcommand{\be}{\begin{equation}} 
\newcommand{\ee}{\end{equation}}
\newcommand{\bea}{\begin{eqnarray}} 
\newcommand{\eea}{\end{eqnarray}}
\newcommand{\bean}{\begin{eqnarray*}} 
\newcommand{\eean}{\end{eqnarray*}}
\newcommand{\rf}[1]{(\ref {#1})}
\newcommand{\un}{{\rm 1\!I}}
\newcommand{\dx}{\,{\rm d}x}
\newcommand{\dy}{\,{\rm d}y}
\newcommand{\ds}{\,{\rm d}s}
\newcommand{\dt}{\,{\rm d}t}
\newcommand{\la}{\langle}
\newcommand{\ra}{\rangle} 

\def\a{\alpha}
\def\g{\gamma}
\def\go{\gamma_0}
\def\uu{u(x,t)\,u(y,t)\,{\rm d}x\,{\rm d}y}
\def\x{\frac{x}{|x|}}
\def\y{\frac{y}{|y|}}
\def\xy{\frac{x-y}{|x-y|}}
\def\axy{\frac{|x|+|y|}{|x-y|}\left(1-\frac{x}{|x|}\cdot\frac{y}{|y|}\right)}
\def\bxy{\left(\frac{x}{|x|}-\frac{y}{|y|}\right)\cdot\frac{x-y}{|x-y|}}

\def\X{\mathcal{X}}
\def\t{\tau}
\def\s{\sigma} 
\def\e{\varepsilon}
\def\p{\partial}
\def\div{\nabla\cdot}
\def\f{\varphi}
\def\r{\varrho}
\def\h{{H^1(\Omega)}}

\title{\Large\textbf{Blowup of solutions \\ to  a diffusive aggregation model}} 
\author{Piotr Biler$^1$,  Grzegorz Karch$^1$, Philippe Lauren\c cot$^2$\\
\small  $^1$ Instytut Matematyczny, Uniwersytet Wroc{\l}awski \\
\small pl. Grunwaldzki 2/4, 50--384 Wroc{\l}aw, Poland\\
\small{\tt \{Piotr.Biler,Grzegorz.Karch\}@math.uni.wroc.pl}\\
\small $^2$ Institut de Math\'ematiques de Toulouse, CNRS UMR~5219 \\
\small Universit\'e de Toulouse, F--31062 Toulouse C\'edex 9, France\\
\small {\tt Philippe.Laurencot@math.univ-toulouse.fr}}

\date{\today}
\maketitle

\begin{abstract} 
\noindent The nonexistence of global in time solutions is studied for a class of 
aggregation equations involving L\'evy diffusion operators and general interaction kernels. 
\end{abstract}

\bigskip

\baselineskip=18pt

\section{Introduction} 
\setcounter{equation}{0}

We consider in this paper the Cauchy problem for the evolution equation 
\bea 
\p_tu+\nu(-\Delta)^{\a/2}u+\div(u\nabla K\ast u)=0,\label{eq}\\
u(x,0)=u_0(x),\label{ini}
\eea
which describes swarming, collective motion and aggregation phenomena in biology and mechanics of continuous media. Here $x\in\R^d$, $t\ge 0$, and $u=u(x,t)\ge 0$ is either the population density of a species or the density of particles in a granular media. 

When $\nu=0$, equation~\rf{eq} can be considered as either a conservation law with a nonlocal (quadratic) nonlinearity or a transport equation with nonlocal velocity, and its character depends strongly on the properties of the kernel $K$. A classical choice for $K$ is $K(x)={\rm e}^{-|x|}$ or, more generally, $K$ is a~radially symmetric function of $r=|x|$. {\em Nonincreasing} kernels correspond to the attraction of particles while {\em nondecreasing} ones are repulsive. Local and global existence of solutions to the inviscid equation \rf{eq}  ($\nu=0$) has been thoroughly studied in \cite{L} under some additional hypotheses on the kernel, see also \cite{BB,BL}. In particular, kernels that are smooth (not  singular) at the origin $x=0$ lead to the global in time existence of solutions, see e.g. \cite{BL,L}. Mildly singular kernels (e.g. $C^1$ off the origin, like $K(x)={\rm e}^{-|x|}$) may lead to blowup of solutions either in finite or infinite time \cite{BB,BCL,BL,L,LR07,LR08}. Strongly singular kernels like potential type (arising in chemotaxis theory, cf. \cite{B-AMSA,BK}) 
\be 
K(x)=c|x|^{\beta-d}, \label{pot}
\ee 
with  $1<\beta<d$  (so, in particular, the Newtonian potential kernel $K(x)=c_d|x|^{2-d}$, $d\ge 3$), usually lead to finite time blowup of ``large'' solutions, see \cite{BK,BW,LRZ08}.

Equations~\rf{eq} with fractional diffusion term ($\nu>0$) have been introduced in the physical literature and studied in, e.g., \cite{BFW,BW,BWu,E}, starting in nineties of the 20th century. The linear term in \rf{eq} is described by a~fractional power of the Laplacian operator in $\R^d$ (or, more generally, by a~L\'evy diffusion operator) defined in the Fourier variables by 
\be
{\mathcal F}((-\Delta)^{\a/2}u)(\xi)= |\xi|^\a\widehat u(\xi),\label{F}
\ee
with $0<\a\le 2$. 
\bigskip

When $\nu>0$ and $K$ is a radially symmetric and nonincreasing function of $r=|x|$ with a mild singularity at $r=0$, equation~\rf{eq} then features a diffusive term which spreads the distribution of particles and a nonlinear drift term which concentrates it, thus acting in the opposite direction. The fundamental question for the Cauchy problem \rf{eq}--\rf{ini} is to decide whether $u$ governed by the competition between the nonlinear transport term and the linear dissipative term, can describe aggregation phenomena or not. Of course, the answer may depend on the (regularity and) size of initial data. 

Typical approaches to prove a finite time aggregation include an extension of the method of {\em characteristics} \cite{BCL,LR08}, the {\em energy} method (e.g. \cite{BB,BL,LR07,LR08}), and the {\em moment} (or {\em virial}) method. The latter has been first applied to mean field models for self-gravitating particles and chemotaxis systems, \cite{B3}, and recently in \cite{B-AMSA,BK}. At this point, we mention that the characteristics method cannot obviously be applied in the presence of diffusion.
\bigskip

Our aim in this paper is to present a simple virial type argument showing finite time blowup of a large class of solutions of \rf{eq}--\rf{ini} if $\nu\ge 0$ and $0<\alpha<1$. The results we obtained are similar to those in \cite{LR08} but we believe that our proofs are more direct and simpler. In addition, our assumptions on the initial data (see \rf{even} and \rf{integrability} below) and the kernel (see \rf{ker0} and \rf{ker} below) are less restrictive than those in \cite{LR07} (radial symmetry and high localization on the initial data, the kernel $K$ being a nonincreasing function of $r=|x|$) and \cite{LR08} (existence of exponential moments or even compact support assumption and $K(x)={\rm e}^{-|x|}$). In particular, the kernel $K$ is allowed to have a repulsive part.

The case of the strong dissipation $1<\a\le 2$ and nonlinearities with potential kernels \rf{pot} has been considered in \cite{BK,BW,LRZ08} where threshold conditions on the values of $\a$, $\beta$, $d$ have been determined so that solutions can be either continued indefinitely in time, or they can blow up in a finite time for suitable initial data.  But, for weakly singular kernels as the ones considered in this paper, the strong dissipation $1<\a\le 2$ prevents finite time blowup and global solutions exist, see \cite{E} and \cite[Theorem~3]{LR07}.

\medskip

\noindent {\bf Notation. } 

\noindent 
The integrals with no integration limits are meant as $\int_{\R^d}\dots$\ . 
Various positive constants are denoted by $C$; sometimes the dependence of $C$ on parameters is written explicitly, e.g. $C=C_\e$, $C=C(\e)$.

\section{Main results}\label{sec:mr}
\setcounter{equation}{0}

We do not consider here local in time existence of solutions, their positivity and mass conservation properties since these topics have been discussed in detail in, e.g., \cite{BB,BK,BW,L}. 
One should note that the existence theory in \cite{LR07, LR08} is developed in the spirit of arguments used for conservation laws in \cite{L}, i.e. without taking into account regularization effects of the diffusion while \cite{BK,BW,BWu,E} employed those effects in a significant way. 
\bigskip

 Our results employ a crucial property of the gradient $\nabla K(x)$ of the convolution kernel in \rf{eq}, namely the fact that $-x/|x|$ is its homogeneous part near the origin. More precisely, we will use two sets of assumptions on the kernel $K$. There is a locally Lipschitz continuous function $k$ such that $K(x)=k(|x|)$ for $x\in\R^d$, and 
\begin{description}
\item[(H1)] either: \ \ there is $K_0>0$ such that 
\be 
-K_0 \le k' \le 0 \;\;\mbox{ and }\;\; \kappa_R:= - \sup_{(0,R)}{ k' } > 0 \label{ker0}
\ee
for each $R>0$; 
\item[(H2)] or: \ \ $k'(r)=-k_1'(r)+k_2'(r)$ for $r>0$ and there are $K_0>0$, $K_1>0$,  $K_2\ge 0$, and $\delta\in [0,1)$ such that 
\begin{equation}
|k'(r)| \le K_0\ (1 + r^\delta)\,, \quad K_1 \le k_1'(r)\,, \quad k_2'(r) \le K_2\ r^\delta \label{ker}
\end{equation}
for all $r>0$.
\end{description}

We will consider solutions which are {\em even} in $x$ which is implied by the assumption that the initial condition $u_0$ is even: 
\be
u_0(x)=u_0(-x), \qquad x\in\R^d, \label{even}
\ee 
together with the radial symmetry of the kernel $K$ and the uniqueness of solutions to \rf{eq}--\rf{ini}. 
Moreover, we need that 
\be 
M_1:= \int_{\R^d} |x|\ u_0(x) \dx < \infty.\label{integrability}
\ee
As we have already remarked, the total mass is conserved during the evolution of \rf{eq}--\rf{ini}  
\be
\int_{\R^d}u(x,t)\dx=M:=\int_{\R^d}u_0(x)\dx.\label{mass}
\ee 

\medskip

Now we are in a position to state main results of this work. 
\medskip 

\begin{thm}\label{th:inviscid} 
Assume that $\nu=0$. Consider a nonnegative and integrable initial condition $u_0\not\equiv 0$ satisfying \rf{even} and \rf{integrability}. 

\begin{description}
\item[(a)]  If $k$ fulfills \rf{ker0}  or \rf{ker} with $K_2=0$, then the solution $u$ to the Cauchy problem \rf{eq}--\rf{ini} ceases to exist in a finite time. 
\item[(b)] If $k$ fulfills \rf{ker} and $M_1$ is sufficiently small, then the  solution $u$ to the Cauchy problem \rf{eq}--\rf{ini} ceases to exist in a finite time. 
\end{description}
\end{thm}

 Let us emphasize here that the assumptions \rf{ker0} and \rf{ker} with $K_2=0$ apply to two different classes of kernels $k$: indeed, $k'$ is required to be bounded in the former but can vanish at infinity (in the sense that $\kappa_R$ might decay to zero as $R\to\infty$). The growth condition is less restrictive for the latter but $k'$ is not allowed to vanish at infinity. We also point out here that we do not know whether the smallness of $M_1$ is a necessary condition for finite time blowup to occur when $k$ fulfills \rf{ker}. 

\medskip

It follows from Theorem~\ref{th:inviscid} that, in the absence of diffusion and if the kernel is attractive ($k'\le 0$), finite time blowup takes place for any nonzero initial data while a partially attractive kernel  seems to require the initial data to be sufficiently concentrated for this phenomenon to occur. Since diffusion is expected to act also as a repulsive term, localization of the initial data seems also to be needed for finite time blowup when $\nu>0$, even if the kernel is attractive. Indeed, we have the following result.

\medskip

\begin{thm}\label{th:dissipative}
Consider a nonnegative and integrable initial condition $u_0\not\equiv 0$ satisfying \rf{even} and \rf{integrability}. Assume that $\nu>0$, $0<\alpha<1$,  and that $k$ fulfills either \rf{ker0} or \rf{ker}. If $M$ is sufficiently large and $M_1$ is sufficiently small, then the solution $u$ to the Cauchy problem \rf{eq}--\rf{ini} ceases to exist in a finite time. 
\end{thm}

\bigskip

Observe that, besides localization of the initial data as in the partially repulsive case, Theorem~\ref{th:dissipative} also requires  the total mass $M$ to be sufficiently large. This is due to the fact that, in the proof, it does not seem to be possible to balance the contribution from the diffusion with that from the drift term.

In contrast to \cite{LR07} and \cite[Theorem~12]{LR08}, our conditions on $u_0$ guaranteeing finite time blowup do not require the $L^1$-norm of $u_0$ to be smaller (in a suitable sense) than that of $u_0 (K * u_0)$. We also improve \cite[Theorem~8]{LR08} where $u_0$ is assumed to be in $L^1(\R^d;e^{2|x|}\dx)$. 

\section{Virial inequalities}\label{sec:vi}
\setcounter{equation}{0}

For $\g\in (0,1]$ and $x\in \R^d$ we define $W_\g(x)=w_\g(|x|)$ with 
$$
w_\g(r)=\frac{1}{\g}\left((1+r)^{\g}-1\right),  \quad r\ge 0.
$$ 
Evidently, $W_\g\ge 0$ is a Lipschitz continuous function which will be used as a weight function. We list below some properties of $W_\g$ and $w_\g$ we will repeatedly use in the remainder of the paper. 

\begin{lem}\label{diss} Consider $\gamma\in (0,1)$ and $\alpha\in (\gamma,1)$. Then $(-\Delta)^{\a/2}W_\g\in L^\infty(\R^d)$.
\end{lem}

\proof Recall that, for $x\in \R^d$, the L\'evy--Khintchine representation formula reads
$$
(-\Delta)^{\a/2}W_\g(x) = - C(d,\alpha)\ \int_{\R^d} \frac{W_\g(x+y)-W_\g(x)}{|y|^{d+\alpha}}\ \dy
$$
with
$$
C(d,\alpha) := \frac{\alpha \Gamma((\alpha+d)/2)}{2 \pi^{(d+2\alpha)/2} \Gamma((2-\alpha)/2)},
$$
see, e.g., \cite[Theorem~1]{DI} 
{or} \cite{J}. Given $x\in \R^d$, we set $\varrho=\varrho(|x|):=\max{\{1,|x|\}}$, and use the monotonicity and subadditivity of $r\mapsto r^\g$ to obtain
\begin{eqnarray*}
\frac{1}{C(d,\alpha)}\ \left| (-\Delta)^{\a/2}W_\g(x) \right| & \le & \frac{1}{\g} \int_{\R^d} \frac{(1+|x|+|y|)^\g - (1+|x|)^\g}{|y|^{d+\a}} \dy \\
& \le & \frac{1}{\g} \int_{B(0,\varrho)} \g (1+|x|)^{\g-1} \ |y|^{1-d-\a} \dy \\
& & + \frac{1}{\g} \int_{\R^d\setminus B(0,\varrho)} \frac{(1+|x|)^\g +|y|^\g - (1+|x|)^\g}{|y|^{d+\a}} \dy \\
& \le & C(d) \ (1+|x|)^{\g-1}\ \frac{\varrho^{1-\a}}{1-\a} + C(d)\ \frac{\varrho^{\g-\a}}{\g (\a-\g)} \\
& \le & C(d,\a,\g)\ \varrho^{\g-\a}\,,
\end{eqnarray*}
and the right-hand side of the above inequality is bounded since $\a>\g$. \qed

\bigskip

Additional properties of $w_\g$ are summarized in the next lemma.

\begin{lem}\label{weight} 
Consider $\g\in (0,1]$. For each $\e>0$ there exists a constant $C_\e>0$ such that the inequalities  
\be 
(1-w_\g'(r))\le (1-\g)\ w_\g(r) \;\;\mbox{ and }\;\; w_\gamma(r) \le \frac{1}{\g} r^\g\le \e+C_\e w_\g(r)\label{w}
\ee
hold for all $r\ge 0$.   For $\delta\in [0,\g)$ and $R>1$ we have 
\be 
r^\delta \le \frac{2 w_\g(r)}{R^{\g-\delta}} \;\;\;\mbox{ for }\;\;\; r\ge R\,.\label{acdc}
\ee 
\end{lem}

\proof 
The first inequality in \rf{w} follows from the observation that the function $f(r)= (1-\g) \left((1+r)^\g-1\right) - \g + \g (1+r)^{\g-1}$ satisfies $f'(r)=\g(1-\g)(1+r)^{\g-2}r\ge 0$ and $f(0)=0$. 

The second inequality in \rf{w} is clear for small $r\ge 0$ and suitably large $C=C_\e$, as well as for large $r\gg 1$. 

 Finally, if $R>1$ and $r\ge R$, we have
$$
w_\g(r) = \int_1^{1+r} s^{\g-1} \ds \ge r\ (1+r)^{\g-1} \ge \frac{(1+r)^\g}{2} \ge \frac{r^{\g-\delta}}{2}\ r^\delta\,,
$$
from which \rf{acdc} readily follows since $\g>\delta$. \qed

\bigskip

Next we derive an identity involving the moment $I_\g$ of a nonnegative solution $u$ of \rf{eq} defined by 
\be
I_\g(t):=\int_{\R^d} W_\g(x) u(x,t) \dx\label{mom}
\ee 
whenever it is meaningful (e.g., if $u\in L^1((0,T);(1+|x|)\dx)$). 

\begin{lem}\label{le:vi} For each $t\ge 0$, we have
\be 
\frac{{\rm d}I_\g}{\dt}(t) = -\nu \mathcal{D}(t) + \mathcal{A}_1(t) + \mathcal{A}_2(t), \label{ev}
\ee
where
\begin{eqnarray*}
\mathcal{D}(t) & := & \int_{\R^d} u(x,t) \left[ (-\Delta)^{\a/2}W_\g\right](x) \dx ,\\
\mathcal{A}_1(t) & := & \iint (w_\g'(|x|)-1) k'(|x-y|)\x\cdot\xy\uu , \\
\mathcal{A}_2(t) & := & \iint k'(|x-y|)\x\cdot\xy\uu .
\end{eqnarray*}
\end{lem}

\proof The evolution of $I_\g$ is governed by \rf{eq} so that 
\bean 
\frac{{\rm d}I_\g}{\dt}(t) & = & -\nu\int_{\R^d} W_\g(x) \left[ (-\Delta)^{\a/2}u \right](x,t)\dx \\
&\quad &+ \int_{\R^d} w_\g'(|x|)\x u(x,t) \cdot (\nabla K\ast u)(x,t)\dx \\
&=& -\nu\int_{\R^d} u(x,t) \left[ (-\Delta)^{\a/2}W_\g\right](x) \dx \\
&\quad&+\iint w_\g'(|x|)k'(|x-y|)\x\cdot\xy\uu \\
&=& -\nu\mathcal{D}(t) \\
&\quad&+\iint (w_\g'(|x|)-1) k'(|x-y|)\x\cdot\xy\uu \\
&\quad&+\iint k'(|x-y|)\x\cdot\xy\uu ,
\eean 
whence \rf{ev}. \qed

\bigskip

The next step is to find suitable upper bounds for $\mathcal{D}$, $\mathcal{A}_1$, and $\mathcal{A}_2$. Such an estimate for $\mathcal{D}$ follows at once from Lemma~\ref{diss} and \rf{mass},  and reads
\bea
\mathcal{D}(t) & \le & \left\| (-\Delta)^{\a/2}W_\g\right\|_\infty \ \int_{\R^d} u(x,t) \dx, \nonumber\\
\mathcal{D}(t) & \le & C(d,\gamma,\alpha) M \;\;\mbox{ \text{provided} }\;\; \alpha\in (\gamma,1) .\label{lin2}
\eea

\medskip

Now we turn to $\mathcal{A}_1$ and $\mathcal{A}_2$, and first consider the case where $k$ satisfies \rf{ker0}.

\begin{lem}\label{le:est1A}
Assume that $k$ fulfills \rf{ker0}. Then, for any $R>1$,
\begin{eqnarray}
\mathcal{A}_1(t) & \le & (1-\gamma)\ M\ K_0\ I_\g(t)\,, \label{e1a1} \\
\mathcal{A}_2(t) & \le & M\ \kappa_{2R}\ \left( \frac{I_\g(t)}{w_\g(R)} - \frac{M}{2} \right)\,. \label{e1a2}
\end{eqnarray}
\end{lem}

\proof We infer from \rf{ker0}, \rf{mass}, and the first inequality in \rf{w} that 
\begin{eqnarray*}
\mathcal{A}_1(t) & \le & \iint (1-w_\g'(|x|))\ |k'(|x-y|)|\ \uu \\
& \le & (1-\g)\ K_0 \iint w_\g(|x|)\ \uu \,,
\end{eqnarray*}
whence \rf{e1a1}. Symmetrizing the double integral $\mathcal{A}_2$, we obtain 
\begin{eqnarray*}
\mathcal{A}_2(t) & = & \frac{1}{2}\ \iint k'(|x-y|)\ \bxy\ \uu \\
& = & \frac{1}{2}\ \iint k'(|x-y|)\ \axy\ \uu \,.
\end{eqnarray*}
Since $k'$ is nonpositive by \rf{ker0} and 
\be
1\le \frac{|x|+|y|}{|x-y|}\,, \qquad (x,y)\in \R^d\times \R^d\,,\label{spirou}
\ee
we deduce from \rf{ker0} and \rf{mass} that, for any $R>1$, 
\begin{eqnarray*}
\mathcal{A}_2(t) & \le & \frac{1}{2}\ \iint k'(|x-y|)\ \left(1-\x\cdot\y\right)\ \uu \\
& \le & \frac{1}{2}\ \int_{B(0,R)} \int_{B(0,R)} k'(|x-y|)\ \left(1-\x\cdot\y\right)\ \uu \\
& \le & - \frac{\kappa_{2R}}{2}\ \int_{B(0,R)} \int_{B(0,R)} \left(1-\x\cdot\y\right)\ \uu \\
& \le & - \frac{\kappa_{2R}}{2}\ \iint \left(1-\x\cdot\y\right)\ \uu \\
& \quad & + \kappa_{2R}\ \int_{\R^d} \!\! \left( \int_{\R^d\setminus B(0,R)} \left(1-\x\cdot\y\right)\ u(y,t) \dy\right) u(x,t) \dx  \\
& \le & -\frac{\kappa_{2R}}{2}\ \left( M^2 - \left( \int \x\ u(x,t)\ \dx \right)^2 \right) \\
& \quad & + \kappa_{2R}\ \int_{\R^d} \!\!\left( \int_{\R^d\setminus B(0,R)} \frac{w_\g(|y|)}{w_\g(R)}\ u(y,t) \dy\right) u(x,t) \dx \\
& \le & -\frac{\kappa_{2R}}{2}\ M^2 + \frac{\kappa_{2R}}{w_\g(R)}\ M\ I_\gamma(t)\,,
\end{eqnarray*}
since the evenness of $x\mapsto u(x,t)$ warrants that 
\be 
\int_{\R^d} \x\ u(x,t)\ \dx = 0 \;\;\mbox{ for }\;\; t\ge 0\,. \label{fantasio}
\ee
The proof of Lemma~\ref{le:est1A} is then complete. \qed

\bigskip

We now derive the counterpart of Lemma~\ref{le:est1A} when $k$ satisfies the weaker assumption \rf{ker}. Though the proof roughly proceeds along the same steps as that of Lemma~\ref{le:est1A}, it is more complicated because some terms involving $k_1$ and $k_2$ have to be handled separately. 
\medskip

\begin{lem}\label{le:est2A}
Assume that $k$ fulfills \rf{ker}. Then, for any $R>1$, $\gamma\in (\delta,1)$ and $\e\in (0,1)$, there is a constant $C_\e'$ depending only on $\e$, $\g$ and $\delta$ such that
\begin{eqnarray}
\mathcal{A}_1(t) & \le &  K_0\ \left[ (1-\g)\ \left( 2M + M\ R^\delta + I_\g(t) \right) + 2\ M\ R^{\delta-\g} \right]\ I_\g(t) \,, \qquad \label{e2a1} \\
\mathcal{A}_2(t) & \le & M \left( \frac{K_1}{w_\g(R)} + K_2\ C_\e' \right) I_\g(t) - \left( \frac{K_1}{2} - 2\ K_2\ \e \right) M^2\,. \label{e2a2}
\end{eqnarray}
\end{lem}

\proof On the one hand, we infer from \rf{ker} and \rf{w} that, for $R>1$, 
\begin{eqnarray*}
\mathcal{A}_1(t) & \le & \iint (1-w_\g'(|x|))\ |k'(|x-y|)|\ \uu \\
& \le & K_0 \iint (1-w_\g'(|x|))\ \left( 1+|x-y|^\delta \right)\ \uu \\
& \le & (1-\g)\ K_0\ \iint w_\g(|x|)\ \left( 1 + |y|^\delta \right)\ \uu \\
&\quad &   + K_0\ \int_{B(0,R)}  (1-w_\g'(|x|))\ |x|^\delta\ u(x,t) \left( \int_{\R^d} u(y,t) \dy \right)\ \dx  \\
&\quad &   + K_0\ \int_{\R^d\setminus B(0,R)} (1-w_\g'(|x|))\ |x|^\delta\ u(x,t) \left( \int_{\R^d} u(y,t) \dy \right)\ \dx  \,.
\end{eqnarray*}
We next use \rf{mass}, \rf{w},   \rf{acdc},  and the property $w_\g'\le 1$ to obtain
\begin{eqnarray*}
\mathcal{A}_1(t) & \le & (1-\g)\ K_0\ \left( \int_{\R^d} \left( 1 + (1+|y|)^\gamma \right)\ u(y,t) \dy \right)\ I_\g(t) \\
&\quad &   + (1-\g)\ M\ K_0\ \int_{B(0,R)} w_\g(|x|)\ |x|^\delta\ u(x,t) \dx \\
&\quad &   + \frac{2\ M\ K_0}{R^{\g-\delta}}\ \int_{\R^d\setminus B(0,R)} w_\g(|x|)\ u(x,t) \dx  \\
& \le & (1-\g)\ K_0\ \left( 2\ M + \g\ I_\g(t) \right)\ I_\g(t) + (1-\g)\ K_0\ R^\delta\ M\ I_\g(t) \\
&\quad &   + \frac{2\ M\ K_0}{R^{\g-\delta}}\ I_\g(t) \\
& \le &   K_0\ \left[ (1-\g)\ \left( 2M + M\ R^\delta + I_\g(t) \right) + 2\ M\ R^{\delta-\g} \right]\ I_\g(t)  \,, 
\end{eqnarray*}
  hence \rf{e2a1}. 
\medskip

On the other hand, after the symmetrization of the double integral in $\mathcal{A}_2$, it follows from the positivity of $k_1'$ in \rf{ker} that, for $R>1$,
\begin{eqnarray*}
\mathcal{A}_2(t) & = & \frac{1}{2}\ \iint k'(|x-y|)\ \bxy\ \uu \\
& = & \frac{1}{2}\ \iint k'(|x-y|)\ \axy\ \uu \\
& = & \frac{1}{2}\ \int_{B(0,R)} \int_{B(0,R)} k'(|x-y|)\ \axy\ \uu \\
&\quad & +\frac{1}{2}\ \int_{B(0,R)} \left( \int_{\R^d\setminus B(0,R)} k'(|x-y|)\ \axy\ u(y,t) \dy\right) u(x,t) \dx \\
&\quad & +\frac{1}{2}\ \int_{\R^d\setminus B(0,R)} \left( \int_{\R^d} k'(|x-y|)\ \axy\ u(y,t) \dy\right) u(x,t) \dx  \\
& \le & - \frac{1}{2}\ \int_{B(0,R)} \int_{B(0,R)} k_1'(|x-y|)\ \axy\ \uu \\
&\quad & +\frac{1}{2}\ \iint k_2'(|x-y|)\ \axy\ \uu \,.
\end{eqnarray*}
Recalling \rf{ker}, \rf{mass}, \rf{spirou} and \rf{fantasio}, we proceed as in the proof of Lemma~\ref{le:est1A} to conclude that
\begin{eqnarray*}
& & \int_{B(0,R)} \int_{B(0,R)} k_1'(|x-y|)\ \axy\ \uu \\
& \ge & K_1\ \int_{B(0,R)} \int_{B(0,R)} \left(1-\x\cdot\y\right)\ \uu \\
& \ge & K_1\ M^2 - 2\ K_1\ \int_{\R^d} \left( \int_{\R^d\setminus B(0,R)} u(y,t) \dy\right) u(x,t) \dx  \\
& \ge & K_1\ M^2 - 2\ \frac{M\ K_1}{w_\g(R)}\ I_\g(t)\,.
\end{eqnarray*}
Using once more \rf{ker}, \rf{mass} and the obvious bound
$$
0\le\axy=\bxy\le 2\,, \qquad (x,y)\in \R^d\times \R^d\,,
$$ 
we find
\begin{eqnarray*}
& & \frac{1}{2}\ \iint k_2'(|x-y|)\ \axy\ \uu \\
& \le & K_2\ \iint |x-y|^\delta\ \uu \\
& \le & K_2\ \iint \left( |x|^\delta+|y|^\delta \right)\ \uu \\
& \le & 2\ M\ K_2\ \int |x|^\delta\ u(x,t) \dx \\
& \le & 2\ M\ K_2\ \left( \e\ M + \frac{\delta}{\g}\ \left( \frac{\g-\delta}{\g \e} \right)^{(\g-\delta)/\delta}\ I_\g(t) \right)\,,
\end{eqnarray*}
whence \rf{e2a2}. \qed

\section{Finite time blowup}\label{sec:bu}
\setcounter{equation}{0}

\subsection{The inviscid case $\nu=0$}\label{sec:bu1}

Now we are ready to prove the first blowup result for the inviscid model \rf{eq}--\rf{ini} with $\nu=0$, and begin with the case of a nonincreasing kernel $K$. 
\bigskip

\noindent{\it Proof of Theorem~\ref{th:inviscid}~(a).\ } 
We argue by contradiction and assume the solution $u$ of \rf{eq}--\rf{ini} to be well-defined for all times. Combining  \rf{ev}, \rf{e1a1} and \rf{e1a2}, we end up with
$$
\frac{{\rm d}I_\g}{\dt}(t) \le \Lambda_{\g,R}(I_\g(t)) := M\ \left [ (1-\gamma)\ K_0 + \frac{\kappa_{2R}}{w_\g(R)} \right]\ I_\g(t) - \frac{\kappa_{2R}\ M^2}{2} 
$$
for all $R>1$ and $\g\in (0,1)$.  Since $\Lambda_{\g,R}$ is a nondecreasing function, we realize that we have $I_\g(t)\le I_\g(0) + \Lambda_{\g,R}(I_\g(0))\ t$ \, for $t\ge 0$ as soon as $\Lambda_{\g,R}(I_\g(0))<0$. Then, of course, $I_\g$ attains zero at some finite time $t_0$ which is impossible for nonnegative regular solutions to \rf{eq}--\rf{ini}, a contradiction  with the  global existence. 

We next observe that we can always find $\g\in (1/2,1)$ and $R>1$ such that $\Lambda_{\g,R}(I_\g(0))<0$ or equivalently
$$
\left [ (1-\gamma)\ K_0 + \frac{\kappa_{2R}}{w_\g(R)} \right]\ I_\g(0) < \frac{\kappa_{2R}\ M}{2}\,. 
$$
Indeed, if $\g\in (1/2,1)$, we have $w_\g(r)\le r$ for $r\ge 0$ and $w_\g(r)\ge \sqrt{r}/2$ for $r\ge 1$. Therefore, choosing $R>1$ such that $M_1<(M\ \sqrt{R})/8$ and then $\g\in (1/2,1)$ such that $(1-\g)\ K_0\ M_1<(\kappa_{2R}\ M)/4$, we realize that 
$$\left [ (1-\gamma)\ K_0 + \frac{\kappa_{2R}}{w_\g(R)} \right]\ I_\g(0) \le \left [ (1-\gamma)\ K_0 + \frac{2\kappa_{2R}}{\sqrt{R}} \right]\ M_1 < \frac{\kappa_{2R}\ M}{2}\,.
$$
With this choice of $R$ and $\gamma$,  we have $\Lambda_{\g,R}(I_\g(0))<0$ and the proof is complete. \qed

\bigskip

\noindent{\it Proof of Theorem~\ref{th:inviscid}~(b).\ } 
We again argue by contradiction and assume the solution $u$ of \rf{eq}--\rf{ini} to be well-defined for all times. Combining  \rf{ev}, \rf{e2a1} and \rf{e2a2}, we end up with
\be
\frac{{\rm d}I_\g}{\dt}(t) \le \Lambda_{\g,R,\e}(I_\g(t)) \;\;\mbox{ for }\;\; t\ge 0\,, \label{spip}
\ee
where
\begin{eqnarray*}
\Lambda_{\g,R,\e}(z) & := & M\ \left[ (1-\g)\ K_0\ \left( 2 + R^\delta \right) +  \frac{2 K_0}{R^{\g-\delta}}   + \frac{K_1}{w_\g(R)} + K_2\ C_\e' \right] z \\
&\quad  & + (1-\g)\ K_0\ z^2 -   \left( \frac{K_1}{2} - 2\ K_2\ \e \right)   M^2
\end{eqnarray*}
for all $R>1$ and $\e\in (0,1)$. As before, the inequality \rf{spip} contradicts the global existence of nonnegative regular solutions to \rf{eq}--\rf{ini} as soon as $\Lambda_{\g,R,\e}(I_\g(0))<0$. Since $\Lambda_{\g,R,\e}$ is an increasing function and $I_\g(0)\le M_1$, we have $\Lambda_{\g,R,\e}(I_\g(0))\le \Lambda_{\g,R,\e}(M_1)$. Observing that an appropriate choice of $\e$ (sufficiently small) and $R$ (sufficiently large) warrants $\Lambda_{\g,R,\e}(0)<0$, we thus have $\Lambda_{\g,R,\e}(M_1)<0$ provided $M_1$ is small enough. Hence, for such a choice of $\e$ and $R$, finite time blowup of the solution to \rf{eq}--\rf{ini} occurs as claimed. 

  Finally, if $K_2=0$, we may argue as at the end of the proof of Theorem~\ref{th:inviscid}~(a) to show that, given any nonzero initial condition $u_0$, we may find $R$ large enough and $\g$ close to one such that $\Lambda_{\g,R,\varepsilon}(I_\g(0))<0$, which completes the proof.~\qed

\bigskip

 Clearly,  the  only term in \rf{spip} that prevents Theorem~\ref{th:inviscid}~(b) from being valid for an arbitrary nonzero initial condition $u_0$ is the term $K_2\ C_\e'\ I_\g$ which cannot be made arbitrarily small by an appropriate choice of $\g$, $R$, and $\e$. This term reflects the deviation of $k$ from being decreasing, and thus the partially repulsive behaviour of $k$.

\subsection{The dissipative case $\nu>0$}\label{sec:bu2}

The second result applies to solutions with suitably large initial data in the dissipative case:
\medskip

\noindent{\it Proof of Theorem~\ref{th:dissipative}.\ } Assume first that $k$ fulfills \rf{ker0}. We argue by contradiction and assume the solution $u$ of \rf{eq}--\rf{ini} to be well-defined for all times $t\ge 0$. Combining  \rf{ev}, \rf{lin2}, \rf{e1a1} and \rf{e1a2}, we end up with
\begin{eqnarray*}
\frac{{\rm d}I_\g}{\dt}(t) & \le & \Lambda_{\g,R}(I_\g(t)) := M\ \left [ (1-\gamma)\ K_0 + \frac{\kappa_{2R}}{w_\g(R)} \right]\ I_\g(t) \\
&\quad & \hspace{3cm}+\ \nu C(d,\g,\alpha)\ M - \frac{\kappa_{2R}\ M^2}{2} 
\end{eqnarray*}
for all $R>1$ and $\g\in (0,\alpha)$.  As before, the above inequality contradicts the global existence of nonnegative regular solutions to \rf{eq}--\rf{ini} as soon as $\Lambda_{\g,R}(I_\g(0))<0$, the latter being true if $\Lambda_{\g,R}(M_1)<0$. Fix $\g\in (0,\alpha)$ and $R>1$ and assume that $M>4\nu\ C(d,\g,\alpha)/\kappa_{2R}$. Then, $\Lambda_{\g,R}(0) \le -M\nu\ C(d,\g,\alpha)<0$ so that $\Lambda_{\g,R}(M_1)<0$ if $M_1$ is sufficiently small.

If $k$ fulfills \rf{ker}, the proof is similar and relies on \rf{ev}, \rf{lin2}, \rf{e2a1} and \rf{e2a2}.\qed

\medskip

In contrast to the proof of Theorem~\ref{th:inviscid}~(a), we cannot play with the parameter $\g$ in the proof of Theorem~\ref{th:dissipative} when $k$ fulfills \rf{ker0}. Indeed, $\g$ is limited by the constraint $\g<\alpha$, and cannot be chosen arbitrarily close to one. This explains the necessity to have sufficiently localized initial data in the sense that $M_1$ is required to be small enough.

\bigskip\bigskip

\noindent {\bf Acknowledgements.}
The preparation of this paper was partially supported by the Polish Ministry of 
Science grant N201 022 32/0902, the POLONIUM projects \'EGIDE no. 13886SG (2008) and (2009),  
and by the European Commission Marie Curie Host Fellowship 
for the Transfer of Knowledge ``Harmonic Analysis, Nonlinear
Analysis and Probability''  MTKD-CT-2004-013389. 
The authors are greatly indebted to Tomasz Cie\'slak for pointing them out the preprint \cite{LR07}. They also thank the referee for comments that improved the earlier version of the paper. 
 
\newpage

\end{document}